\magnification=\magstep1
\input amstex
\documentstyle{amsppt}

\NoBlackBoxes
\loadbold
\def\mbox#1{\leavevmode\hbox{#1}}


\def\diag{\operatorname{diag}}
\def\Domain{\operatorname{Domain}}

\def\dvol{\operatorname{dvol}}
\def\End{\operatorname{End}}
\def\ev{\operatorname{ev}}
\def\id{\operatorname{id}}
\def\Index{\operatorname{Index}}
\def\Ker{\operatorname{Ker}}

\def\opp{\operatorname{op}}

\def\Range{\operatorname{Range}}
\def\supp{\operatorname{supp}}
\def\spec{\operatorname{spec}}

\def\gtimes{\hat{\otimes}}
\def\bsigma{\boldsymbol\sigma}
\def\E{\Cal{E}}
\def\b{\Bbb}

\def\R{\Bbb R}
\def\C{\Bbb C}
\def\H{\Cal{H}}

\def\KK{\Cal{K}}
\def\L{\Cal{L}}

\def\lcl{[\![}
\def\rcl{]\!]}


\topmatter
\title On Graded $K$-theory, Elliptic Operators and the Functional Calculus \endtitle
\rightheadtext{Graded $K$-theory, Functional Calculus, and Elliptic Operators}
\author Jody Trout \endauthor
\thanks Research of the author is partially supported by NSF grant DMS-9706767\endthanks
\address 
Department of Mathematics \newline
Dartmouth College \newline
Hanover, NH 03755  \endaddress
\email jody.trout\@dartmouth.edu \endemail

\abstract  Let $A$ be a graded $C^*$-algebra. We characterize Kasparov's $K$-theory group
$\hat{K}_0(A)$ in terms of graded $*$-homomorphisms by proving a general converse to the functional calculus theorem for
self-adjoint regular operators on graded Hilbert modules. An application to the index theory of elliptic differential operators on smooth
closed manifolds and asymptotic morphisms is discussed.
\endabstract

\endtopmatter


\document

\head {\bf 1. Introduction.} \endhead

Let $A$ be a graded $\sigma$-unital $C^*$-algebra with grading automorphism $\alpha$. We characterize Kasparov's $K$-group in the 
category of graded $C^*$-algebras, $\hat K_0(A)=KK(\C,A)$, as the group of  graded homotopy classes
of graded $*$-homomorphisms from $C_0(\R)$,  the $C^*$-algebra of continuous functions on the real line with the even-odd function
grading,  to the graded tensor product $A \gtimes \KK$, where $\KK \cong M_2(\KK)$ is the $C^*$-algebra of compact operators graded
into diagonal and off-diagonal matrices. Addition is given by direct sum.

The isomorphism is established in Section 3 by proving a general converse to the functional calculus theorem \cite{11} for self-adjoint
regular operators on graded Hilbert modules in Section 2. We will indicate in Section 4 how this characterization is useful in simplifying
calculations with asymptotic morphisms of $C^*$-algebras and elliptic differential operators $D$ with coefficients in a trivially graded
$C^*$-algebra $A$ over a smooth closed manifold $M$. The functional calculus will give an explicit formulation as (nontrivial)
compatible graded $*$-homomorphisms of the generalized Fredholm index  $\Index_A(D) \in K_0(A)$ and the symbol class $[\sigma(D)]
\in K^0_A(T^*M)$ (the topological $K$-theory of vector $A$-bundles of the cotangent bundle $T^*M$) in a form
which is suitable for composing directly with asymptotic morphisms,  with no rescaling or suspensions as in the general theory. Since the
product in $E$-theory is given by composition, this approach to index theory is simpler than using the Kasparov product in $KK$-theory
\cite{10}, which can be very technical. 

We should note that Kasparov's graded $K$-theory is unrelated to van Daele's version, except when $A$ is trivially graded \cite{19}. This
paper represents work that partially began in the author's thesis \cite{17}, although the material in Section 2 is new. The author would
like to thank his advisers Nigel Higson and Paul Baum for their invaluable help and encouragement and also Erik Guentner for helpful
suggestions.

\head 2. {\bf Graded $C^*$-algebras and Hilbert modules.} \endhead

In this section we collect some definitions and results on graded $C^*$-algebras and Hilbert modules and fix notation. For a complete
discussion, see the books \cite{3,9}.

Let $A$ be a $C^*$-algebra. Recall that $A$ is {\it graded} if there is a $*$-automorphism $\alpha : A \to A$ such that 
$\alpha^2 = \id_A$.  Equivalently, there is a decomposition of $A$ as a direct sum $A = A_0 \oplus A_1,$ where
$A_0$ and $A_1$ are self-adjoint closed linear subspaces with the property that if $x \in A_m$ and $y \in A_n$ then $xy \in A_{m + n} \
\pmod{2}$. In fact,  $A_n = \{x \in A : \alpha(x) = (-1)^n x\}.$  We write $\partial x = n$ if $x \in A_n$. If there is a self-adjoint unitary
$\epsilon$ (called the {\it grading operator}) in the multiplier algebra $M(A)$ such that
$\alpha(x) = \epsilon x \epsilon^*$, then $A$ is said to be {\it evenly graded}. A $*$-homomorphism $\phi : A \to B$ of graded
$C^*$-algebras is  {\it graded} if $\phi(A_n) \subset B_n$ for $n = 0,1$.

\example{Example 2.1} The following are the main examples we will be concerned with.
\roster
\item"a.)" Every $C^*$-algebra $A$ can be {\it trivially} graded by setting
$A_0 = A$ and $A_1 = \{0\}$. This is an even grading with grading operator
$\epsilon = 1$. The complex numbers $\C$ are always assumed to be trivially graded.
\item"b.)" The $C^*$-algebra $C_0(\R)$ of continuous complex-valued functions on $\R$ vanishing at infinity is graded into the
even and odd functions by defining
$\alpha(f)(t) = f(-t)$ for all functions $f \in C_0(\R)$ and $t \in \R$.
\item"c.)" Let $\H$ be an infinite-dimensional separable Hilbert space. By choosing an isomorphism $\H \cong \H \oplus \H$ we obtain
the standard even grading on the
$C^*$-algebra of compact operators $\KK = \KK(\H) \cong M_2(\KK)$, with grading operator $$\epsilon = \pmatrix 1 & 0 \\ 0 & -1
\endpmatrix$$ which is determined uniquely up to conjugation by a unitary homotopic to the identity.
\endroster
\endexample

Let $A$ and $B$ be graded $C^*$-algebras. Define a graded product and graded involution on the vector space tensor product $A \odot B$
by the formulas
$$\aligned (a \gtimes b)(a' \gtimes b') &= (-1)^{\partial b \partial a'}(aa'
\gtimes bb') \\ (a \gtimes b)^* &= (-1)^{\partial a \partial b}(a^* \gtimes b^*).
\endaligned$$  The resulting $*$-algebra is denoted $A \hat{\odot} B$. A grading on $A \hat{\odot} B$ is defined by setting
$$\partial(a \gtimes b) = \partial a + \partial b \ \pmod{2}.$$  Now faithfully represent $A$ and $B$ by $\rho_1$ and $\rho_2$ on
graded Hilbert spaces $H_1$ and $H_2$ with grading operators
$\epsilon_1$ and $\epsilon_2$, respectively. Then $A \hat{\odot} B$ is faithfully represented on $H_1 \otimes H_2$ (graded by
$\epsilon_1 \otimes \epsilon_2$) via the formula 
$$\rho(a \gtimes b) = \rho_1(a)\epsilon_1^{\partial a} \otimes  \rho_2(b).$$ The $C^*$-algebra completion is denoted $A \gtimes B$ and is
called the (minimal) {\it graded} tensor product. It does not depend on the choice of representations. (There is also a maximal graded
tensor product \cite{3} but it will not be needed for our purposes since one of the factors will always be nuclear.)

\proclaim{Lemma 2.2}{\rm (Proposition 15.5.1 \cite{3})} If $B$ is evenly graded, then $A \gtimes B \cong A \otimes B$. If
$A$ is also evenly graded, then under this isomorphism $A \otimes B$ is also evenly graded.
\endproclaim

\proclaim{Corollary 2.3} Let $\KK$ have the standard even grading.  Then $A \gtimes \KK \cong M_2(A \otimes \KK)$. If $A$ is evenly
graded by $\epsilon$, $A \gtimes \KK
\cong M_2(A \otimes \KK)$ with standard even grading given by 
$\eta = \diag (\epsilon \otimes 1 ,  -\epsilon \otimes 1 ).$
\endproclaim

Let $B$ be another graded $C^*$-algebra with grading $\beta$. Then $B[0,1] = C([0,1], B)$ canonically inherits a grading by the formula
$\hat{\beta}(f)(t) = \beta(f(t))$. Two graded $*$-homomorphisms $\phi_0, \phi_1 : A \to B$ are {\it graded homotopic} if there is a graded
$*$-homomorphism $\Phi : A \to B[0,1]$ such that composition with the evaluation maps $\ev_t : B[0,1] \to B$ for $t = 0, 1$ are equal to
$\phi_0$ and $\phi_1$, respectively.  We shall denote by $\lcl A, B\rcl$ the set of {\it graded} homotopy classes of {\it graded}
$*$-homomorphisms from $A$ to $B$. If $\phi : A \to B$ is a graded
$*$-homomorphism, then we denote by $\lcl \phi \rcl$ its equivalence class in $\lcl A, B \rcl$.

A Hilbert $A$-module $\H$ is {\it graded} if there is a Banach space decomposition $\H = \H_0 \oplus \H_1$ such that $\H_n \cdot A_m
\subseteq \H_{n+m}$  and $\langle \H_n, \H_m \rangle \subseteq A_{n+m}$ $\pmod{2}$. We let $\L(\H)$ denote the $C^*$-algebra of all
bounded $A$-linear maps $T : \H \to \H$ with an adjoint $T^*$ and let
$\KK(\H)$ denote the closed two-sided ideal of compact operators. The grading on $\H$ induces gradings on $\L(\H)$ and $\KK(\H)$ via
the identities $\partial T = m$ if $T(\H_n) \subset \H_{n+m}$. We let $\H^{\opp}$ denote $\H$ with the opposite grading $\H^{\opp}_n =
\H_{1-n}$. Note that if $A$ is trivially graded, $\H$ is the direct sum of two orthogonal $A$-modules. If $\phi : B
\to \L(\H)$ is a $*$-homomorphism,  a closed submodule $\E$ of $\H$ is {\it $\phi$-invariant} if $\phi(b) : \E \to \E$ for all $b \in  B$.

\head {\bf 3. The Converse Functional Calculus.} \endhead

Let $\H$ be a (graded) Hilbert $A$-module.  A {\it regular operator} on $\H$ is a densely defined closed $A$-linear map $D : \Domain(D)
\to \H$ such that the adjoint $D^*$ is densely defined and $1 + D^*D$ has dense range.  \footnote{If $\H = A$ then $D$ is sometimes called
an unbounded multiplier \cite{4,8,11}.} $D$ has {\it degree one} if $\partial(Dx) = \partial x + 1$ for all $x \in \Domain(D)$.  

\proclaim{Proposition 3.1} For any graded $*$-homomorphism $\phi : C_0(\R) \to A$, there is a maximal $\phi$-invariant closed graded
Hilbert $A$-submodule $A_\phi$ of $A$ and a self-adjoint regular operator $D$ on $A_\phi$ of degree one such that for all $f \in C_0(\R)$
we have $\phi(f)|_{A_\phi} = f(D).$
\endproclaim

\demo{Proof} Given a graded $*$-homomorphism $\phi : C_0(\R) \to A$, define $$A_\phi = C_0(\R) \gtimes_\phi A =
\overline{\phi(C_0(\R))A}$$ to be the closed right ideal generated by the image of $\phi$. This is a closed graded Hilbert submodule of
$A$ (see Blackadar \cite{3}.) Let $C_c(\R)$ denote the dense graded ideal of continuous functions on $\R$ with {\it compact support}.
Define $$\Domain(D) = \phi(C_c(\R))A$$ which is a dense graded submodule of $A_\phi$. Let $d$ denote the function $d(t) = t$ on $\R$.
Define $D : \Domain(D) \to A_\phi$ by the formula
$D\phi(f)x = \phi(df)x$ where $f \in C_c(\R)$ (so $df \in C_c(\R)$) and extend linearly. Suppose that
$\phi(f)x = \phi(g)y$ for some other $g \in C_c(\R)$. Choose a function $d'
\in C_c(\R)$ such that $d = d'$ on the compact set $\supp(f) \cup \supp(g)$. Then we have
$$D\phi(f)x = \phi(d'f)x = \phi(d')\phi(f)x = \phi(d')\phi(g)y = \phi(d'g)y = D\phi(g)y.$$  It follows that $D$ is well-defined and is clearly
$A$-linear. Also, $D$ is degree one since $d$ is an odd function on $\R$. The computation 
$$\langle D\phi(f)x, \phi(g)y\rangle = x^*\phi(df)^*\phi(g)y = x^*\phi(d\bar{f}g) = x^*\phi(f)^*\phi(dg)y = \langle \phi(f)x, D\phi(g)y
\rangle$$ shows that $D$ is symmetric on $\Domain(D)$. This implies that $D$ is closeable, so we replace $D$ by its closure $\bar{D}$.
Consequently, $(D \pm i)$ are injective and have closed range by Lemma 9.7 \cite{11}. Let $f \in C_c(\R)$. For any $x \in A$ we have
$$(1 + D^2)\phi((1 + d^2)^{-1})\phi(f)x = \phi((1+d^2)(1+d^2)^{-1}f)x = \phi(f)x.$$ 
It follows that $\Range(1 + D^2) \supseteq \Domain(D)$ is dense and so $D$ is regular. We will show $D$ is self-adjoint by using a Cayley
transform argument.

Extend $\phi$ to $\phi^+ : C_0(\R)^+ \to A^+$ by adjoining a unit. Let $z \in C_0(\R)^+$ denote the unitary 
$$z(t) = \frac{t - i}{t+i} = 1 - 2ir_-(t), \text{ for } t \in \R$$
where $r_-(t) = (t -i)^{-1}$ denotes the resolvent. Let $U_D = \phi^+(z) = 1 - 2i \phi(r_-) \in A^+$. It is easy to check that
for all $x \in \Domain(D)$ we have that the unitary $U_D$ satisfies
$$U_D(D + i)x = (D + i)U_Dx = (D - i)x.$$
By Lemma 9.8 and the discussion following Proposition 10.6 in Lance \cite{11}, the closed symmetric regular operator $D$ is
self-adjoint and $U_D = (D+i)^{-1}(D-i)$).

To show $\phi(f)|_{A_\phi} = f(D)$, it suffices to show this for the resolvents $r_\pm(t) = (d \pm i)^{-1}(t)$. Let
$\{f_n\}_{n=1}^{\infty}$ be an approximate unit for
$C_0(\R)$ consisting of compactly supported functions.  Let $x \in A_\phi$ be given. Then $\phi(f_n)x \in
\Domain(D)$ for all $n$ and $\phi(f_n)x  \to x$ as $n \to \infty$. We have that as $n \to \infty$,
$$(D \pm i)\phi((d \pm i)^{-1}f_n)x = \phi((d\pm i)(d \pm i)^{-1}f_n)x = \phi(f_n)x \to x.$$  Now since
$\phi((d \pm i)^{-1}f_n)x = \phi((d \pm i)^{-1})\phi(f_n)x \to \phi((d \pm i)^{-1})x$
 as $n \to \infty$ and $(D \pm i)$ is closed, we conclude that $\phi((d \pm i)^{-1})x = (D \pm i)^{-1}x$. Since $x \in A_\phi$ was arbitrary,
we are done. \hfill \qed 
\enddemo

Let $B$ be a $C^*$-algebra. If $\H$ is a Hilbert $B$-module, a $*$-homomorphism $\phi : A \to \L(\H)$ is called {\it
nondegenerate} if $\phi(A)\H$ is dense in $\H$. It is called {\it strict} if $\{\phi(u_n)\}$ is Cauchy in the strict topology of $\L(\H)$ for
some approximate unit $\{u_n\}$ in $A$. Nondegeneracy implies strictness \cite{11}.  The following result may be considered the
converse to the functional calculus for self-adjoint regular operators \cite{2,4,11}.

\proclaim{Theorem 3.2} {\rm (Converse Functional Calculus)} Let $\phi : C_0(\R) \to \L(\H)$ be graded. There is a closed graded
$\phi$-invariant Hilbert submodule $\H_\phi$ of $\H$ and a self-adjoint regular operator $D$ on $\H_\phi$ of degree one such that for all
$f \in C_0(\R)$ we have
$\phi(f)x = f(D)x$ for all $x \in \H_\phi$. Moreover, if $\phi$ is strict then $\H_\phi$ is complemented and $\phi(f) = f(D) \in \L(\H_\phi)
\subseteq \L(\H)$. If $\phi$ is nondegenerate then $\H = \H_\phi$.  And if $\phi(C_0(\R)) \subset \KK(\H)$ then $D$ has compact
resolvents.
\endproclaim

\demo{Proof} Let $A = \L(\H)$. Let $D' : \Domain(D') \to A_\phi$ be the self-adjoint regular operator on $A_\phi = C_0(\R) \gtimes_\phi
A$ from the previous proposition such that $\phi(f) = f(D')$.  Let $i : A \to \L(\H)$ be the identity.  Define $\H_\phi =
\overline{\phi(C_0(\R)\H} = A_\phi \gtimes_i \H$ which is a  closed Hilbert submodule of $\H$.  Define $D = D' \gtimes_i 1$ on 
$$\Domain(D) = \Domain(D') \hat{\odot}_i \H \supseteq \phi(C_c(\R))\H.$$  By Proposition 10.7 \cite{11}, $D$ extends to a self-adjoint
regular operator on $\H_\phi$.  ($D = i_*(D')$ in the notation of \cite{11}.) If $x \in \H_\phi$, we compute that
$$f(D)x = f(D' \gtimes_i 1)x = (f(D')
\gtimes_i 1)x = f(D') \gtimes_{i} x  = \phi(f)x.$$ If $\phi$ is strict then $\H_\phi$ is a complemented submodule of $\H$ by
Proposition 5.8
\cite{11} and so $\L(\H_\phi)$ includes as a graded subalgebra of $\L(\H)$.  The result now easily follows.
\qed \enddemo

Note that  if $\phi$ is the zero homomorphism then $\H_\phi = \{0\}$ and $D = 0$, so $f(D) = 0 = \phi(f)$.

\head {\bf 4. Graded $K$-theory.} \endhead

\demo{\bf Standing Assumptions} Throughout this section $A$ will denote a complex $\sigma$-unital graded $C^*$-algebra and
$C_0(\R)$ and $\KK$ will have the gradings as in Example 2.1. \enddemo

Let $H_A$ denote the Hilbert $A$-module of all sequences $\{a_n\}_1^\infty \subset A$ such
that $\{\sum_1^n a_k^*a_k\}_{1}^\infty$ converges in $A$. It has a natural grading into sequences of even and odd elements. Let
$\hat{H}_A = H_A \oplus H^{\opp}_A$, where
$H^{\opp}_A$ denotes $H_A$ with the {\it  opposite} grading. This is the  standard graded Hilbert module for $A$. We have the following
very important result of Kasparov in the theory of graded Hilbert modules.

\proclaim{Stabilization Theorem} {\rm (Kasparov \cite{10})} If $\H$ is a countably generated graded Hilbert $A$-module then $\H
\oplus \hat{H}_A \cong \hat{H}_A$. \endproclaim

It is a standard result that $A \gtimes \KK$ is graded $*$-isomorphic to $\KK(\hat{H}_A)$, the $C^*$-algebra of compact operators on
$\hat{H}_A$ (with the induced grading) (See 14.7.1 \cite{3}). For the remainder of this section, we will identify $A \gtimes \KK =
\KK(\hat{H}_A)$. From stabilization, conjugation by the graded isomorphism $\hat{H}_A \cong \hat{H}_A \oplus \hat{H}_A$ determines a
unitary in
$\L(\hat{H}_A) = M(A \gtimes \KK)$ of degree zero.

\proclaim{Lemma 4.1} Let $u \in M(A \gtimes \KK)$ be a unitary of degree zero. There is a strictly continuous path of degree zero
unitaries $\{U_t\}_{t \in [0,1]} \subset M(A \gtimes \KK)$ such that $U_1 = u$ and $U_0 = 1$. \endproclaim

\demo{Proof} Write $\KK = \KK(H \oplus H)$ where $H = L^2[0,1]$. Then $M(A \gtimes \KK)$ contains a copy of $\L(H \oplus H)$. Let
$\{v_t\}$ be a strictly continuous path of isometries in $\L(H)$ with $p_t = v_tv_t^* \to 0$ strongly as $t \to 0$ as in Proposition 12.2.2
\cite{3} . Set $V_t = v_t \oplus v_t \in \L(H \oplus H)$ and note that each $V_t$ has degree zero. Set $W_t = 1 \gtimes V_t$ which also
has degree zero and let 
$$U_t = W_t u W_t^* + (1 - W_t W_t^*)$$
for $t > 0$ and $U_0 = 1$. It is easy to check that this works. \hfill \qed \enddemo

\demo{\bf Definition 4.2} Let $A$ have grading automorphism $\alpha$. Define 
$$K'(A) = K'(A, \alpha) = \lcl C_0(\R), A \gtimes \KK \rcl.$$ 
Define a binary operation on $K'(A)$ by direct sum $\lcl \phi \rcl + \lcl \psi \rcl = \lcl \phi \oplus \psi \rcl$, where the direct sum is with
respect to the graded isomorphism $\hat{H}_A \cong \hat{H}_A \oplus \hat{H}_A$ \enddemo

\proclaim{Theorem 4.3} $K'(A)$ is an abelian group under the direct sum operation and satisfies the relation
$$-\lcl\phi\rcl = \lcl u \phi u^* \rcl$$ where $u = u^* = \pmatrix 0 & 1 \\ 1 & 0
\endpmatrix$ on $\hat{H}_A = H_A \oplus H_A$.
\endproclaim

\demo{Proof}  It follows from Lemma 4.1 and the proof of Theorem 3.1 in Rosenberg \cite{15} carried over to the graded case that
$K'(A)$ is an abelian monoid with zero given by the zero (or any null-homotopic) $*$-homomorphism. We only need to show inverses 

Let $\phi : C_0(\R) \to \KK(\hat{H}_A)$ be graded. Let $D$ be the regular operator on $\H_\phi \subset \hat{H}_A$ associated to $
\phi $ from the converse functional calculus. Via stabilization $\H_\phi \oplus \hat{H}_A \cong \hat{H}_A$ and Lemma 4.1, we may
assume (up to graded homotopy) that $\phi$ is strict by Proposition 5.8 \cite{11}. Thus $\phi(f) = f(D)$ for all $f \in C_0(\R)$.  Then
$D^{\opp} = uDu^*$ on the Hilbert module
$\H_\phi$ is the operator associated to $\lcl u\phi u^* \rcl$ since by the functional calculus
$$f(D^{\opp}) = f(uDu^*) = uf(D)u^* = u\phi(f)u^*.$$  Let $\epsilon$ be the grading on $\hat{H}_A$. For each $t \geq 0$, define
 $$\b D_t = \pmatrix D & t \epsilon \\ t
\epsilon & D^{\opp} \endpmatrix$$ on $\H_\phi \oplus \H_\phi^{\opp} \subseteq \hat{H}_A$ and let $\b D_t = 0$ on the
 complement. Define $\Phi_t : C_0(\R) \to \KK(\hat{H}_A)$ by
$$\Phi_t(f) = f(\b D_t).$$ For $t = 0$ we have $\Phi_0(f) = f(\b D_0) = \phi \oplus u \phi u^*$. Note that
$$\b D_t^2 = \pmatrix D & t\epsilon \\ t\epsilon & D^{\opp} \endpmatrix^2 =
\pmatrix D^2 + t^2 & 0 \\ 0 & D^{\opp}{}^2 + t^2 \endpmatrix$$ and so the spectrum of $\b D_t$ is contained outside the interval $(-t ,t)$.
Therefore, 
$$\|f(\b D_t)\| \leq \sup \{|f(x)| : x \in \spec(\b D_t)\} \to 0 \text{ as } t \to \infty$$ 
for all $f \in C_0(\R)$ and the result follows. \hfill \qed \enddemo 

\demo{\bf Definition 4.4} A {\it $K$-cycle} for a graded $C^*$-algebra $A$ is an ordered pair $(\H, T)$, such that $\H = \H_0 \oplus \H_1$
is a countably generated graded Hilbert $A$-module and $T \in \L(\H)$, where $\L(\H)$ is the graded $C^*$-algebra of all bounded
$A$-linear operators on $\H$ with adjoint, which satisfies the following conditions:
\roster
\item"i.)" $T$ is of degree one;
\item"ii.)" $T - T^* \in \KK(\H)$ is compact;
\item"iii.)" $T^2 - 1 \in \KK(\H)$ is compact.
\endroster
\noindent The $K$-cycle is called {\it degenerate} if $T^2 = 1$.
\enddemo
By a standard argument we may assume that $T = T^*$ is self-adjoint. There is an obvious notion of {\it unitary equivalence} for two
$K$-cycles \cite{3,10}. Two $K$-cycles $(\H_0, T_0)$ and $(\H_1, T_1)$ are {\it homotopic} if there is a $K$-cycle $(\H, T)$ for
$A[0,1]$ such that
$(\H \gtimes_{\ev_i} A, T \gtimes_{\ev_i} 1)$ are unitarily equivalent to $(\H_i, T_i)$ where $\ev_t : A[0,1] \to A$ are the evaluation
maps. A collection $\{(\H, T_t)\}_{t \in [0,1]}$ of $K$-cycles for $A$ is called an {\it operator homotopy} if $t \mapsto T_t$
is norm continuous in $t$. An operator homotopy induces a homotopy $(\H', T)$ by defining $\H' = C([0,1], \H)$ and $T(f)(t) = T_t(f(t))$ for
$f : [0,1] \to \H$.

\proclaim{\bf Proposition 4.5} {\rm (Theorem 4.1 \cite{10})} The set $KK(\C, A)$ of all equivalence classes
of
$K$-cycles for
$A$ under the equivalence relation (generated by) homotopy is an abelian group under the relations 
$$\aligned (\H_1, T_1) + (\H_2, T_2) &= (\H_1 \oplus \H_2, D_1 \oplus D_2) \\ -(\H,T) &= (\H^{\opp}, -T). \endaligned$$ The class of
any degenerate
$K$-cycle is zero in $KK(\C,A)$.
\endproclaim

Let  $u = \pmatrix 0 & 1 \\ 1 & 0 \endpmatrix$ be the degree one unitary with respect to the grading on $\H = \H_0 \oplus \H_1$

\proclaim{Lemma 4.6} $-(\H, T) = (\H, T^{\opp}) \in KK(\C,A)$, where $T^{\opp} = uTu^*$. \endproclaim

\demo{Proof} In the complex world,  $(\H, T) = (\H, -T)$ since they are operator homotopic (but not through self-adjoint
$K$-cycles in general.) It follows that
$$-(\H,T) = (\H^{\opp},-T) = (\H^{\opp},T) = (\H,uTu^*) = (\H,T^{\opp})$$
since $u : \H^{\opp} \to \H$ implements a unitary equivalence. \hfill \qed \enddemo

\proclaim{Theorem 4.7} $K'(A)$ is isomorphic to $KK(\C, A)$. \endproclaim

\demo{Proof} Let $G(t) = t(t^2 + 1)^{-1/2}$ which defines a degree one, self-adjoint element in $C_b(\R) = M(C_0(\R))$, the continuous
{\it bounded} functions on $\R$. Define a map
$K'(A) \to KK(\C, A)$ via
$$\lcl \phi \rcl \mapsto (\H_\phi, G(D))$$ where $D$ is the regular operator associated to $\phi : C_0(\R) \to 
\KK(\H_\phi) \subset \KK(\hat{H}_A)$ via the converse functional calculus. (As in Theorem 4.3, we may assume that $\phi$ is strict.)
The operator
$G(D)$ is a degree one, self-adjoint element of $M(\KK(\hat{H}_A)) = \L(\hat{H}_A)$ and $G(D)^2 - 1$ is compact since
$$G(D)^2 - 1 = (D^2 + 1)^{-1} = \phi(G) \in \KK(\H_\phi).$$ This map is easily seen to be well-defined since applying the construction to a
graded homotopy $\Phi : C_0(\R) \to \KK(\hat{H}_A)[0,1]$ yields a homotopy of $K$-cycles by using the graded isomorphism 
$$\KK(\hat{H}_A)[0,1] \cong (A \gtimes \KK)[0,1] \cong A[0,1] \gtimes \KK \cong \KK(\hat{H}_{A[0,1]}).$$ It is also distributes over direct
sums and maps
$$-\lcl\phi\rcl = \lcl u \phi u^* \rcl \mapsto (\H_\phi, G(D)^{\opp}) = -(\H_\phi, G(D))$$ via properties of the functional calculus and
Lemma 4.6.

The reverse map is defined using the techniques of Baaj and Julg~\cite{2}. Let
$(\H, F)$ be a $K$-cycle for $A$. We may assume that $F = F^*$ and $\H = \hat{H}_A$. Let $T > 0$  be a strictly positive element of
$\KK(\hat{H}_A)$ of degree zero which commutes with $F$. Any two such operators are operator homotopic via the straight line
homotopy. Let $D = F T^{-1}$. Note that
$\Domain(D) = \Range(T)$ is a dense submodule of $\hat{H}_A$. One has that $D = D^*$ and
$(D^2 + 1)^{-1} = T^2(F^2 + T^2)^{-1}$ is compact. We have the identity
$G(D) = F(F^2 + T^2)^{-1/2}$ and so it also follows that $(\hat{H}_A, F)$ and
$(\hat{H}_A, G(D))$ are operator homotopic. It follows from the identity
$$(D \pm i)^{-1} =  D(D^2 + 1)^{-1} \mp i(D^2 + 1)^{-1}$$ that the resolvents are also compact. Define $$KK(\C, A) \to K'(A)$$ by sending
$(\hat{H}_A, F)$ to the graded homotopy class of the graded $*$-homomorphism
$$\phi : f \mapsto f(D) \in \KK(\hat{H_A}).$$ As above,  $\KK(\hat{H}_{A[0,1]}) \cong \KK(\hat{H}_A)[0,1]$, so a homotopy
$(\hat{H}_{A[0,1]}, F)$ is mapped to a homotopy $\Phi : C_0(\R) \to \KK(\hat{H}_A)[0,1]$. Thus the reverse map is well-defined. One
checks easily that these two maps are inverses to each other.
\qed \enddemo

If $A$ is trivially graded and unital then $A \gtimes \KK  \cong M_2(A \otimes \KK)$  with even grading given by $\epsilon = \diag(1,
-1)$. That is, $M_2(A \otimes \KK)$ is graded into diagonal and off-diagonal matrices. It follows from the above that $$K'(A) = \lcl C_0(\R),
A \otimes \KK\rcl \cong K_0(A).$$ 

We will describe the isomorphism directly via the more familiar language of projections. It is a
standard result that $K_0(A)$ is the group of formal differences of homotopy classes of projections $p = p^* = p^2 \in A \otimes \KK$ with
addition given by direct sum $[p] + [q] = [p' + q']$ where $p
\sim_h p', q \sim_h q'$ and $p' \perp q'$. Let $u \in M_2(M(A \otimes \KK))$ be the degree one unitary 
$$u = \pmatrix 0 & 1 \\ 1 & 0 \endpmatrix.$$ Recall that for any self-adjoint involution $w$ (i.e.,
$w^* = w, w^2 = 1$) there is an associated projection $p(w) = \frac{1}{2}(w + 1).$

Let $x = [p] - [q] \in K_0(A)$ where $p$ and $q$ are projections in
$A \otimes \KK$. Define a map
$$\phi_x : C_0(\R) \to M_2(A \otimes \KK)$$ by the formula
$$\phi_x(f) = \pmatrix f(0)p & 0 \\ 0 & f(0)q \endpmatrix, \quad f \in C_0(\R).$$ This defines a
$*$-homomorphism since
$p = p^2 = p^*$ (similarly for $q$) and is graded since $f(0) = 0$ for any odd function. Note that the homotopy class of $\phi_x$ depends
only on the homotopy classes of $p$ and $q$. Now we define a map $\mu : K_0(A) \to K'(A)$ by mapping
$$x \longmapsto \lcl \phi_x \rcl.$$ It also follows that
$$\phi_{[p]}(f) \oplus \phi_{[q]}(f) = \pmatrix f(0)\diag(p,q) & \diag(0,0) \\ \diag(0,0) & \diag(0,0) \endpmatrix \sim_h \pmatrix f(0)(p' +
q') & 0
\\ 0 & 0
\endpmatrix = \phi_{[p' + q']}(f)$$ and so it is additive. For $x = [p] - [q]$ we have 
$-x = [q] - [p]$ maps to $$\phi_{-x}(f) = \pmatrix f(0)q & 0 \\ 0 & f(0)p \endpmatrix = u \phi_x(f) u^*.$$ Thus, $\mu(-x) = \lcl u\phi_x u^*
\rcl = - \lcl \phi_x \rcl = - \mu(x)$. One should note that with the grading present $\phi_x$ and
$\phi_{-x}$ are {\it not} homotopic through {\it graded} $*$-homomorphisms since $u$ has degree one and the identity has degree zero. 

Conversely, given $\lcl \phi \rcl \in K'(A)$, extend $\phi$ to a graded $*$-homomorphism
$$\phi^+ : C_0(\R)^+ \to (A \otimes \KK)^+$$ by adjoining a unit. Let $z$ denote the unitary given by the ``Cayley transform''
$$z(t) = \frac{t + i}{t-i} = 1 + 2 i r_-(t)$$ where
$r_-(t) = (t - i)^{-1}$ is the resolvent function. Let $u_\phi$ denote the unitary
$$u_\phi = \phi^+(z) = 1 + 2 i \phi(r_-) \in (A \otimes \KK)^+$$   A simple computation shows that 
$(\epsilon u_\phi)^2 = 1 \text{ and } (\epsilon u_\phi)^* = \epsilon u_\phi.$ We also have that $\epsilon^* = \epsilon$ and $\epsilon^2 = 1.$
Consider the associated projections $$p(\epsilon), \ p(\epsilon u_\phi) \in (A \otimes \KK)^+.$$ By the definition of $u_\phi$ above, we see
that $p(\epsilon) - p(\epsilon u_\phi) = 2 i \phi(r_+) \in A \otimes \KK.$ Also, a homotopy of $\phi$ induces a homotopy of the unitary
$u_\phi$ and thus of $p(\epsilon u_\phi)$. We define $\nu : K'(A) \to K_0(A)$ by
$$\nu(\lcl \phi \rcl) = [p(\epsilon)] - [p(\epsilon u_\phi)] \in K_0(A).$$  A simple computation shows that $\nu \circ \mu = 1$. We only
need to show $\mu$ is onto. It then follows that $\nu = \mu^{-1}$ is a homomorphism.

Since $A$ is trivially graded $\hat{H}_A = H_A \oplus H_A$ with each factor determining the grading. Again identify $A \gtimes \KK =
\KK(\hat{H}_A)$. Let $\lcl \phi \rcl \in K'(A)$. Up to graded homotopy we may assume that
$\phi : C_0(\R) \to \KK(\hat{H}_A)$ is strict (via stabilization). Let $$D = \pmatrix 0 & D_+^* \\ D_+ & 0 \endpmatrix$$ on $\hat{H}_A$ be
the self-adjoint regular operator of degree one with compact resolvents from the converse functional calculus such that $\phi(f) = f(D)$. Let
$G(D) = D(D^2 + 1)^{-\frac12}$ which is a self-adjoint bounded operator of degree one on
$\hat{H_A}$ with $G(D)^2 - 1$ compact. By a graded homotopy, we may assume  that $\phi(f) = (f\circ G)(D) = f(G(D))$. (Note that the
diffeomorphism $G : \R \to (-1,1)$ is the homotopy inverse to the inclusion $(-1,1) \subset \R$.) Thus, we can write
$$G(D) = \pmatrix 0 & G_+^* \\ G_+ & 0 \endpmatrix$$  on $H_A \oplus H_A$ where $G_+ : H_A \to H_A$ is a generalized
Fredholm operator \cite{18}. Up to a compact perturbation of $G_+$ (which would induce a graded homotopy), we may assume that
$\Ker(G(D)) = \Ker(G_+) \oplus \Ker(G_+^*)$ is a finite projective $A$-module in $\hat{H}_A$, and is thus complemented. Note that for $x
\in \Ker(G(D))$ we have $f(G(D))x = f(0)x$.  Since $A$ is trivially
graded, $\Ker(G_+)$ and $\Ker(G_+^*)$ are finite projective $A$-modules. Let $P_+^{(*)} \in \KK(H_A)$ be the compact
projections onto
$\Ker(G_+^{(*)})$. Let
$x = [P_+] - [P_+^*] =
\Index_A(G_+)  \in K_0(A)$ \cite{18}. A graded homotopy connecting $\phi$ to the graded $*$-homomorphism
$$\phi_x(f) = \pmatrix f(0)P+ & 0 \\ 0 & f(0)P_+^* \endpmatrix \in \KK(H_A \oplus H_A) = \KK(\hat{H}_A)$$ is given by
$$\Phi_t(f) = \cases f(t^{-1}G(D)), &t > 0 \\ \phi_x(f), &t = 0. \endcases$$
Thus, $\mu (x) = \lcl \phi \rcl$ and so $\mu$ is onto as was desired.

\proclaim{Corollary 4.8} If $A$ is unital and trivially graded then the maps $\mu$ and $\nu$ are inverses. 
\endproclaim


\head {\bf 5. Elliptic Operators over $C^*$-algebras.} \endhead

In this section, we will show how the previous results and the functional calculus give explicit realizatons as graded $*$-homomorphisms
of the $K$-theory symbol class and Fredholm index of an elliptic differential operator with coefficients in a trivially graded $C^*$-algebra.

Let $A$ be a trivially graded {\it unital} $C^*$-algebra and $M$ be a smooth closed Riemannian manifold. Let $E \to M$ and $F \to M$ be
smooth vector $A$-bundles, that is, smooth locally trivial fiber bundles on $M$ whose fibers $E_p$ and $F_p$ are finite projective
$A$-modules for each $p \in M$.  Let $C^\infty(E)$ denote the vector space of smooth sections of $E$, which is a module over $A$,
similarly for $C^\infty(F)$.  Let  $D : C^\infty(E) \to C^\infty(F)$ be an elliptic differential $A$-operator of order $n$ on $M$ \cite{13,17}. (If
$A = \C$ then $D$ is an ordinary differential operator.) Let $\sigma = \sigma(D) : \pi^*(E) \to \pi^*(F)$ denote the principal symbol of $D$
which is a homomorphism of vector
$A$-bundles, where $\pi : T^*M \to M$ is the cotangent bundle. The condition of ellipticity is the requirement that for each
non-zero cotangent vector $\xi \neq 0 \in T_p^*M$ the principal symbol $\sigma_\xi(D) : E_p \to F_p$ is an isomorphism of $A$-modules.

Equipping the fibers $E_p$ (and $F_p$) with smoothly varying Hilbert $A$-module structures $$\langle \cdot , \cdot \rangle_p : E_p
\times E_p\to A$$ defines a pre-Hilbert $A$-module structure on $C^\infty(E)$ via the formula
$$\langle s, s'\rangle = \int_M \langle s(p), s'(p) \rangle_p\dvol_M \in A,$$
for $s, s' \in C^\infty(E),$ where $\dvol_M$ is the Riemannian volume measure on $M$.  
(And any two such structures are homotopic via the straight line homotopy.) It follows that an adjoint differential operator $D^t :
C^\infty(F) \to C^\infty(E)$ exists and is of the same order as
$D$. The principal symbol of the adjoint is  the adjoint of the principal symbol $\sigma_\xi(D^t) = \sigma^*_\xi(D) \in \L(F_p,E_p)$ for $\xi
\in T_p^*M$. Consider the formally self-adjoint differential
$A$-operator of degree one
$$\b{D} = \pmatrix 0 & D^t \\ D & 0 \endpmatrix : C^\infty(E) \oplus C^\infty(F) \to C^\infty(E) \oplus C^\infty(F)$$ 
on the graded pre-Hilbert $A$-module $C^\infty(E) \oplus C^\infty(F)$.
The principal symbol of $\b D$ is the self-adjoint bundle morphism of degree one
$$\bsigma = \bsigma(\b D) = \pmatrix 0 & \sigma^* \\ \sigma & 0 \endpmatrix : \pi^*(E)
\oplus\pi^*(F) \to\pi^*(E) \oplus \pi^*(F)$$ on the graded pull-back vector $A$-bundle $\pi^*(E) \oplus \pi^*(F)$.

\proclaim{Lemma 5.1} The resolvents $(\bsigma \pm i)^{-1} : \pi^*(E) \oplus \pi^*(F) \to \pi^*(E) \oplus \pi^*(F)$ are vector $A$-bundle
morphisms which vanish at infinity on
$T^*M$ in the operator norm induced by the Hilbert $A$-module structures on the fibers $E_p \oplus F_p$.
\endproclaim

\demo{Proof} Follows from homogeneity $\bsigma(p,t\xi) = t^n \bsigma(p,\xi)$ and ellipticity. \qed \enddemo

Form the Cayley transform \cite{14}
$$\bold u = (\bsigma + i)(\bsigma - i)^{-1} = 1 + 2i(\bsigma - i)^{-1}.$$ 
By complementing the vector $A$-bundles $E$ and $F$ , e.g. $E \oplus G \cong M \times A^n$, we may embed
$\pi^*(E \oplus F)$ in a trivial $A$-bundle 
$$\b A = T^*M \times (A^{n}\oplus A^{n}).$$  Now extend the automorphism $\bold u$ to the $A$-bundle $\b A$ by defining it to be equal
to the identity on the complement of $\pi^* E \oplus \pi^* F$ in $\b A$. From the lemma above, it follows that $\bold u$ extends
continuously to the trivial $A$-bundle on the one-point compactification $(T^*M)^+$ by setting $\bold u(\infty) = I$.

Let $\bold \epsilon = \diag (1 , -1 )$ be the grading of the trivial $A$-bundle $(T^*M)^+ \times (A^{n} \oplus
A^{n})$. Since $\bold \epsilon \bsigma = - \bsigma \bold \epsilon$ it follows, as in the previous section,  
$(\bold u \bold \epsilon)^2 = 1$ and $(\bold u \bold \epsilon)^* = \bold u \bold \epsilon$. (We also have obviously that
$\epsilon^* = \epsilon$ and $\epsilon^2 = 1$.)

Therefore, we obtain two projection-valued sections 
$$p(\bold \epsilon) , \ p(\bold u \bold \epsilon) : (T^*M)^+ \to \End(\b A)$$ on $(T^*M)^+$ which are equal at infinity. We can view
them as projection-valued functions $(T^*M)^+ \to M_2(M_n(A)) \cong M_{2n}(A)$. Both define elements in
$K_0(C(T^*M^+) \otimes A)$ and so their difference defines an element
$$\Sigma(D) = [p(\bold \epsilon)] - [p(\bold \epsilon \bold u)] \in K_0(C_0(T^*M) \otimes A).$$ This is the {\it symbol class}
of the elliptic $A$-operator $D$ as constructed in \cite{7,14,17}. By Corollary  4.8 and stability, it follows
that 
$$K_0(C_0(T^*M) \otimes A) \cong \lcl C_0(\R), C_0(T^*M) \otimes M_{2n}(A)) \rcl$$ and $\Sigma(D)$ is
identified with the graded homotopy class of the graded $*$-homomorphism $$\Phi_\sigma : C_0(\R) \to C_0(T^*M, M_{2n}(A)) \cong
M_{2n}(C_0(T^*M) \otimes A)$$ given fiber-wise by the ordinary matrix functional calculus
$$\Phi_\sigma(f)(\xi) = f(\bsigma_\xi(\b D)) \in M_{2n}(A)), \text{ for } \xi \in T^*M.$$

The principal symbol $\sigma(D) : \pi^*(E) \to \pi^*(F)$ determines a class $[\sigma(D)] \in K^0_A(T^*M)$ (the topological
$K$-theory of $T^*M$ defined via vector $A$-bundles) since it is a bundle morphism that is an isomorphism off the compact zero-section
$M \subset T^*M$. By the Mingo-Serre-Swan Theorem \cite{12,16}, we have
$K^0_A(T^*M) \cong K_0(C_0(T^*M) \otimes A)$, which is induced via the action of taking sections as for the case $A = \C$. It thus follows
from this and the constructions in the previous section that all three of these symbol classes can be identified.

\proclaim{Proposition 5.2} $[\sigma(D)] = \Sigma(D) =  \lcl \Phi_\sigma \rcl \in K^0_A(T^*M) \cong K_0(C_0(T^*M) \otimes A).$
\endproclaim

Let $L^2(E)$ denote the completion of the pre-Hilbert $A$-module $C^\infty(E)$. The differential
$A$-operator $\b D$ defines an (essentially) self-adjoint regular operator of degree one on the graded Hilbert $A$-module $\H_D = L^2(E)
\oplus L^2(F)$. (We replace $\b D$ by its closure $\bar{\b D}$ which is self-adjoint.) Since $\b D$ is elliptic, the resolvents $(\b
D \pm i)^{-1}$ are compact. (This follows from the parallel Sobolev theory for differential $A$-operators \cite{13}.) The
complementation of the bundles $E$ and $F$ above allows a coherent inclusion (with the previous constructions)
$$\H_D \subset L^2(\b A) \cong L^2(M) \otimes A^{2n}$$
which induces a graded inclusion of $C^*$-algebras $\KK(\H_D) \hookrightarrow M_{2n}(\KK \otimes A)$.  By the
functional calculus for self-adjoint regular operators
\cite{11} we obtain a graded $*$-homomorphism 
$$\Phi_D : C_0(\R) \to M_{2n}(\KK \otimes A)  :  f \mapsto f(\b D)$$
Recall that the usual definition of the generalized Fredholm (analytic) index $\Index_A(D)$ in terms of kernel and cokernel modules
requires compact perturbations for a general $C^*$-algebra $A$ \cite{13,18}. This is incorporated in
the computations in the proof of Corollary 4.8, so we see that the functional calculus for $\b D$ gives this index.

\proclaim{Proposition 5.3} $\Index_A(D) = \lcl \Phi_D \rcl \in K_0(A)$. \endproclaim

Naturally associated to $M$ and $A$ is an asymptotic morphism of $C^*$-algebras
$$\{\Psi_t\}_{t \in [1,\infty)} : C_0(T^*M) \otimes A \to \KK(L^2M) \otimes A,$$ which is defined via Fourier transforms and a partition of
unity up to asymptotic equivalence. 
(For complete details on
the construction see
\cite{5,7,17}.)  The induced map
$$\Psi_*: K^0_A(T^*M) \cong K_0(C_0(T^*M) \otimes A) \to K_0(A)$$
on $K$-theory is useful for doing index-theoretic and
$K$-theoretic calculations with elliptic operators. If $M = \R^n$,
the induced map is Bott periodicity $K_0(C_0(\R^{2n}) \otimes A) \cong K_0(A)$
\cite{17}. The following result implies the exact form of the Mishchenko-Fomenko index theorem \cite{13}, hence the Atiyah-Singer
index theorem \cite{1} when $A = \C$ as proved originally by Higson \cite{7}.

\proclaim{Theorem  5.4}{\rm (Lemma 4.6 \cite{17})} If $D$ is an elliptic differential $A$-operator of order one on $M$ then
$$\Psi_*([\sigma(D)]) = \Index_A(D) \in K_0(A).$$ \endproclaim

The proof reduces to composing the graded symbol homomorphism
$$\Phi_\sigma : C_0(\R) \to M_{2n}(C_0(T^*M) \otimes A) : f \mapsto f(\bsigma) $$
with the matrix inflation of this ``fundamental'' asymptotic morphism for $M$ and $A$
$$\{\Psi_t\} : M_{2n}(C_0(T^*M) \otimes A) \to M_{2n}(\KK \otimes A).$$
and comparing this to the  continuous family of graded operator $*$-homomorphisms
$$\{\Phi^t_D\}_{t \in [1,\infty)} : C_0(\R)  \to M_{2n}(A \otimes \KK) : f \mapsto f(t^{-1}\b D). $$
One then proves \cite{17} via Fourier analysis and a compactness argument that for any $f \in C_0(\R)$,
$$\lim_{t \to \infty} \| \Psi_t(f(\bsigma)) - f(t^{-1}\b D) \| = 0$$
and so the composition $\{\Psi_t \circ \Phi_\sigma\}$ is asymptotically equivalent to $\{\Phi_D^t\}$. Therefore, by stability and
homotopy invariance of the induced map \cite{5,6}, 
$$\Psi_* \lcl \Phi_\sigma \rcl =  \lcl \Phi^t_D \rcl = \lcl \Phi_D \rcl \in  K_0(A).$$ 
The result now follows by Propositions 5.2 and 5.3.

\bigskip


\enddocument
\end